June 2017

# Explanation of a Polynomial Identity

Nicholas Phat Nguyen[1]

**Abstract.** In this note, we provide a conceptual explanation of a well-known polynomial identity used in algebraic number theory.

A basic theorem in algebraic number theory is that in a number field $E = \mathbb{Q}(w)$ of degree n, with minimal polynomial f(X) for $w$, the dual lattice of $\mathbb{Z}[w]$ relative to the trace form is $\mathbb{Z}[w]/f'(w)$, where $f'(w)$ is the derivative of f(X) evaluated at $X = w$. From that we can deduce that the different ideal of the ring $\mathbb{Z}[w]$ is the principal ideal generated by the number $f'(w)$.

All the standard accounts rely on the following polynomial identity credited to Euler:

$\sum f(X)/(X-u)f'(u) = 1$, where the sum runs over all the conjugates $u$ of $w$,

or more generally,

$\sum q(u)f(X)/(X-u)f'(u) = q(X)$ for any polynomial q of degree $< n = \deg(f)$.

See, e.g., [1], chapter 3. The above polynomial identities are beautiful and can be proven simply by observing that both sides are polynomials of degree $< n = \deg(f)$ and have the same value when X is set equal to any of the n different conjugates of w. However, the identities seem to fall out of the sky and therefore rather mysterious. A conceptual explanation along the following line may be illuminating.

---

[1] E-mail address: nicholas.pn@gmail.com

We start with $E = \mathbb{Q}(w)$, which we can identify with $\mathbb{Q}[X]/f(X)$ via the natural isomorphism that maps X mod f(X) to *w*. E is an algebra over $\mathbb{Q}$, and if we extend the domain of rationality to the algebraic closure K of $\mathbb{Q}$, that is to say if we look at the K-algebra $E \otimes K$ (tensoring over $\mathbb{Q}$), then that K- algebra is isomorphic to a product of n copies of K because over K the polynomial f(X) is a product of n distinct linear factors. The homomorphism from $E \otimes K = K[X]/f(X)$ to a factor K is simply induced by $X \mapsto$ a conjugate of *w* in K.

So we can identify $E \otimes K$ with $\prod K(u)$, where the product runs over the conjugates *u* of *w*, and the symbol $K(u)$ for a conjugate element *u* of *w* refers to the field K obtained by the homomorphism from K[X] to K given by $X \mapsto u$.

The K-algebra $\prod K(u)$ has a natural basis $(\mathbf{e}_u)$, where each vector $\mathbf{e}_u$ has *u*-component 1 and all other components being zero. The sum of all these vectors is equal to the unit element of the K-algebra $\prod K(u)$, namely the vector all of whose components are 1.

What are the elements in $E \otimes K = K[X]/f(X)$ that correspond to this nice natural basis $(\mathbf{e}_u)$ of $\prod K(u)$? A polynomial $\mathbf{p}_u$ in K[X] that maps to $\mathbf{e}_u$ must be divisible by $(X - v)$ for any conjugate *v* of *w* that is $\neq u$. That means $\mathbf{p}_u$ must be divisible by the product of all those factors $(X - v)$ with $v \neq u$, which is just $f(X)/(X - u)$. Moreover, such a polynomial $\mathbf{p}_u$ must leave a remainder of 1 upon division by $(X - u)$, i.e., $\mathbf{p}_u(u) = 1$. The polynomial $f(X)/(X - u)$ when evaluated at *u* has the value $f'(u)$, so the polynomial $f(X)/(X - u)f'(u)$ will map to $\mathbf{e}_u$. Any two such polynomials are congruent mod f(X), so $f(X)/(X - u)f'(u)$ is in fact the only polynomial with degree < deg(f) that maps to $\mathbf{e}_u$. We will denote this polynomial as $\mathbf{p}_u(X)$.

Accordingly, the Euler's polynomial identity $\Sigma f(X)/(X - u)f'(u) = 1$ is an expression of the fact that the sum of $\mathbf{p}_u(X)$ maps to the sum of $\mathbf{e}_u$, and so must be congruent to 1 mod f(X).

The monomial X in K[X] maps to the vector in $\prod K(u)$ whose components are the conjugates of *w*. It follows that for any polynomial q(X) with coefficients in the base field $\mathbb{Q}$, the element q(X) mod f(X) in K[X]/f(X), which corresponds to $q(w) \otimes 1$ in the algebra

E⊗K, maps to the vector whose components are the conjugates of $q(w)$. On the other hand, the sum $\Sigma q(u)f(X)/(X - u)f'(u) = \Sigma q(u) p_u(X)$ also maps to $\Sigma q(u) e_u$ which is the vector whose components are the conjugates of $q(w)$. So $q(X)$ and $\Sigma f(X)q(u)/(X - u)f'(u)$ must be congruent modulo $f(X)$ in the ring $K[X]$. If $q(X)$ has degree $<\deg(f)$, then the expressions must be equal.

______________________________________________________________.

**REFERENCES**:

[1] Serge Lang, *Algebraic Number Theory* (Graduate Texts in Mathematics 110), Springer-Verlag (1986).